\documentclass[12pt]{amsart}
\usepackage{amsmath}
\usepackage{amssymb}
\usepackage{epsfig}
\usepackage{graphicx}
\usepackage{color}
\definecolor{shadecolor}{gray}{0.875}
\usepackage{amscd}
\usepackage{comment}
\usepackage{enumitem}
\usepackage{mathrsfs}

\usepackage[colorlinks=true]{hyperref}

\numberwithin{equation}{section}

\input xy
\xyoption{all}

\calclayout
\allowdisplaybreaks[3]

\theoremstyle{plain}
\newtheorem{prop}{Proposition}[section]

\newtheorem{theo}[prop]{Theorem}

\newtheorem{lemm}[prop]{Lemma}

\theoremstyle{definition}
\newtheorem{defi}[prop]{Definition}

\newtheorem{rema}[prop]{Remark}

\newtheorem{prin}[prop]{Principle}

\def\rH{{\mathrm H}}

\def\Eff{\overline{\mathrm{Eff}}}
\def\Pic{\mathrm{Pic}}

\def\Aut{\mathrm{Aut}}
\def\Nef{\mathrm{Nef}}

\def\Supp{\mathrm{Supp}}

\makeatother
\makeatletter

\author{Brian Lehmann}
\address{Department of Mathematics \\
Boston College  \\
Chestnut Hill, MA \, \, 02467}
\email{lehmannb@bc.edu}

\author{Sho Tanimoto}
\address{Graduate School of Mathematics, Nagoya University, Furocho Chikusa-ku, Nagoya, 464-8602, Japan}
\email{sho.tanimoto@math.nagoya-u.ac.jp}

\title[Manin's conjecture for almost Fano varieties]{A counterexample to Manin's Conjecture for almost Fano varieties}

\begin{document}
\date{\today}

\begin{abstract}
We give an example of a smooth rational threefold over a number field with a big anticanonical divisor that does not satisfy the thin set version of Manin's Conjecture.  In particular, Manin's Conjecture does not hold for rationally connected almost Fano varieties. 
\end{abstract}

\maketitle


\section{Introduction}

Manin's Conjecture predicts the asymptotic growth rate of rational points with bounded anticanonical height on certain projective varieties over global fields.  The precise formulation of this conjecture is a topic of active research.  The statement has been refined several times by various mathematicians including \cite{FMT89}, \cite{BM}, \cite{Peyre}, \cite{BT}, \cite{Peyre03}, \cite{Peyre17}, \cite{LST18}, and \cite{LS24}; in particular, \cite{Peyre03} and \cite{LST18} recognized that one may need to remove a thin set of points in order to obtain the expected growth rate.

Our goal is to analyze the range of applicability of Manin's Conjecture: which types of projective variety can satisfy the desired asymptotic formula?  Classically one considers smooth Fano varieties equipped with the anticanonical height.  There is an alternative first proposed by \cite{Peyre03}:

\begin{defi}
    An almost Fano variety $X$ is a smooth projective geometrically irreducible variety such that $-K_{X}$ is big, $\rH^{1}(X,\mathcal{O}_{X}) = \rH^{2}(X,\mathcal{O}_{X}) = 0$, and $\Pic(X_{\overline{k}})$ is torsion-free.
\end{defi}

Almost Fano varieties have a number of desirable arithmetic properties; in particular, one can define Peyre's constant for such varieties.  This class has the additional advantage that it is broad enough to include the class of toric varieties (\cite{BT98, CLT01}) and equivariant compactifications of vector groups (\cite{CLT02}) as well as many other interesting examples.  On the other hand, \cite[Section 5.1]{LST18} showed that almost Fano varieties can have certain pathological geometric features with surprising arithmetic consequences.

Our main theorem produces a missing counterexample: the thin set version of Manin's Conjecture does not hold for rationally connected almost Fano varieties. (The conjecture remains open for smooth Fano varieties.)  Manin's Conjecture predicts that after removing a thin set the number of rational points of anticanonical height at most $B$ grows like $CB(\log B)^{\rho(X)-1}$.  In our example the actual growth rate is strictly larger than the predicted growth rate no matter what thin set of points we remove.

\begin{theo} \label{theo:maintheorem}
    Let $k$ be any number field.
    There is an example of a smooth rational threefold $X$ over $k$ with $-K_{X}$ big such that for every thin subset $Z \subset X(k)$ there are constants $\mathsf C>0$ and $\epsilon > 0$ satisfying
    \begin{equation*}
        \#\{ \, x \in X(k) \backslash Z \, \,  | \, \, \mathsf H_{-K_{X}}(x) \leq B \, \} > \mathsf CB^{1+\epsilon}
    \end{equation*}
    for every sufficiently large $B$.
\end{theo}

The example $X$ is obtained as a $\mathbb{P}^{1}$-bundle over a rational surface $S$.  This surface $S$ has special geometry: it admits infinitely many $(-1)$-curves and infinitely many elliptic fibrations.  It is constructed as a quotient of a K3 surface using the techniques of \cite{Oguiso89,DO19,DOY23}.  By leveraging the work of \cite{CDJLZ22}, we show that the set of rational points lying on the $(-1)$-curves in $S$ is not thin.

We then construct a $\mathbb{P}^{1}$-bundle $X$ over $S$ with the additional property that $-K_{X}$ is big.  Note that $X$ still admits a countable union of ``negatively embedded'' subvarieties $T_{i}$ which collectively contribute a non-thin set of rational points.  Each $T_{i}$ is a Hirzebruch surface; we verify that the actual growth rate of rational points on each $T_{i}$ is higher than the predicted rate for $X$.  Our construction guarantees that no thin set can discount the contributions from every $T_{i}$ simultaneously.

\bigskip

\noindent
{\bf Acknowledgements:}
The authors would like to thank Dan Loughran and Tim Santens for helpful discussions regarding equidistribution results for non-anticanonical heights. The authors would like to thank Shigeyuki Kondo for discussions regarding Kummer surfaces.
The authors also would like to thank Keiji Oguiso for his inspiring lecture series at Nagoya University in Spring 2022, from which the second author learned the examples in \cite{DO19, DOY23}.
This work was initiated during the authors' visit to American Institute of Mathematics to attend SQuaRE workshop. The authors would like to thank AIM for its hospitality.
The authors would like to thank Runxuan Gao, Dan Loughran, Brendan Hassett, and Tim Santens for comments on our draft.

Brian Lehmann was supported by Simons Foundation grant Award Number 851129.  Sho Tanimoto was partially supported by JST FOREST program Grant number JPMJFR212Z and by JSPS KAKENHI Grant-in-Aid (B) 23K25764.

\bigskip

\noindent
{\bf AI disclosure:}
At a dinner in Pasadena, the authors tried to come up with projects suitable for AI use.  We remembered our old question from \cite{LST18} and discussed the problem over dinner.  At the end, we realized that we already had a solid plan with no need for additional input.

We used AI to look up references.  We also used AI to check the correctness of the paper and AI pointed out a few minor issues. Those are fixed by the authors. None of the ideas, arguments, or writing is due to AI.

\section{Preliminaries}

Throughout we will work over a ground field $k$.  A variety is an integral separated scheme of finite type over the ground field.  Given a locally free sheaf $\mathcal{E}$ on a variety $X$, the notation $\mathbb{P}_{X}(\mathcal{E})$ denotes the projective bundle of quotients of $\mathcal{E}$.

\subsection{Curves and divisors}
Suppose $X$ is a smooth projective variety.  We define $N^{1}(X)$ to be the $\mathbb{R}$-vector space of Cartier divisors up to numerical equivalence.  We will consider the following cones in $N^1(X)$:
\begin{itemize}
    \item $\Eff^1(X)$ is the closure of the cone generated by effective Cartier divisors.
    \item $\Nef^1(X)$ is the nef cone, i.e.~the cone of all $\mathbb R$-Cartier divisors $D$ such that $D.C\geq 0$ for every curve $C$ on $X$.
\end{itemize}
These cones are closed, convex, full-dimensional, and pointed.

\subsection{Thin sets}

We will need the notion of a thin set introduced by Serre in the context of the inverse Galois problem.

\begin{defi}
    Let $X$ be a projective variety over a field $k$.  A thin map is a morphism of projective varieties $f: Y \to X$ such that
    \begin{enumerate}
        \item $f$ is generically finite onto its image, and
        \item $f$ is not birational.
    \end{enumerate}    
    A thin set in $X(k)$ is any subset of a finite union $\cup_{i=1}^{r} f_{i}(Y_{i}(k))$ where $\{f_{i}: Y_{i} \to X\}$ is a finite collection of thin maps.  A thin set coming from a single map $f: Y \to X$ is called type I if $f$ is non-dominant or type II if $f$ is dominant and non-birational.
\end{defi}

This notion is studied in \cite{Serre} and has been extensively studied recently due to its role in the (weak) Hilbert property.

\subsection{Equidistribution and thin sets}

Emmanuel Peyre was the first to consider an equidistribution property in the context of Manin's conjecture (\cite{Peyre}).  This property has been further studied by various mathematicians, e.g., \cite{BL19} and \cite[Section 9.4]{LS24}.  While equidistribution has been mainly considered for anticanonical heights, we will need a formulation for non-anticanonical heights.

Let us recall the setup. Let $X$ be a smooth projective variety defined over a number field $k$ and $L$ be a big divisor on $X$ with an adelic metrization. One can define the Fujita invariant $a(X, L)$ as in \cite[Definition 4.1.1]{Tani26}:
\begin{equation*}
a(X,L) = \inf \{ t \in \mathbb{R} | K_{X} + tL \in \Eff^{1}(X) \}.    
\end{equation*}
We say $(X, L)$ is adjoint rigid if $a(X, L)L + K_X$ has Iitaka dimension $0$.
Let us assume this property. In particular $a(X, L)L + K_X$ is $\mathbb Q$-linearly equivalent to a unique effective $\mathbb Q$-divisor $E$.  We set $X^\circ = X \setminus \mathrm{Supp}(E)$.

Let us fix an adelic metrization for $K_X$ as well. This induces the Tamagawa measure
\[
\tau_{X, E} = \mathsf H_E(x)^{-1}\tau_{X^\circ}
\]
on $X^\circ(\mathbb A_k)$. (See \cite{CLT10} for the definition of this measure.)

The following principle is due to Emmanuel Peyre for anticanonical heights (\cite{Peyre}) and Daniel Loughran and Tim Santens for non-anticanonical heights (\cite{LS24}). Recall that a continuity set is a measurable set whose boundary has measure zero.

\begin{prin}[Peyre, Loughran--Santens]
\label{conj:equidistribution}
   Assume that the geometric fundamental group $\pi_1^{\text{\'et}}(X^\circ)$ is trivial and $\mathrm{Br}(X^\circ)/\mathrm{Br}(k)$ is finite.  Then there exists a thin set $\mathsf Z$ such that for any continuity set $W \subset X^\circ(\mathbb A_k)$, we have
   \[
   \lim_{B \to \infty}\frac{\#\{ x \in X(k) \setminus \mathsf Z \, | \, x \in W, \, \, \mathsf H_L(x) \leq B\}}{\#\{ x \in X(k) \setminus \mathsf Z \, | \,  \mathsf H_L(x) \leq B\}}
   = \frac{\tau_{X, E}(W\cap X^\circ (\mathbb A_k)^{\mathrm{Br}(X^\circ)})}{\tau_{X, E}(X^\circ (\mathbb A_k)^{\mathrm{Br}(X^\circ)})}
   \]
\end{prin}

We will need the following result of Santens establishing equidistribution for toric varieties.

\begin{theo}[{\cite[Theorem 8.11 and Remark 8.12]{Santens26}}]
\label{theo:santens}
    Assume that $X$ is a smooth toric projective variety. Then Principle~\ref{conj:equidistribution} holds with $\mathsf Z$ as the set of rational points on the boundary.
\end{theo}
\begin{proof}
    This follows from \cite[Theorem 8.11 and Remark 8.12]{Santens26} as well as the Portmanteau Theorem, e.g., \cite[Theorem 13.16]{Klenke}.
\end{proof}

Finally we will need the following generalization of \cite[Theorem 1.2]{BL19} and \cite[Theorem 9.13]{LS24}:

\begin{theo}
\label{theo:nocontributionfromthin}
    Assume that $\mathrm{Br}(X^\circ)/\mathrm{Br}(k)$ is finite and Principle~\ref{conj:equidistribution} holds. Let $Z \subset X(k)$ be a thin set on $X$ such that $\mathsf Z \subset Z$.
    Then we have
    \[
\lim_{B \to \infty}\frac{\#\{ x \in Z \setminus \mathsf Z \, | \,  \mathsf H_L(x) \leq B\}}{\#\{ x \in X(k) \setminus \mathsf Z \, | \,  \mathsf H_L(x) \leq B\}} = 0.
\]
   
\end{theo}

\begin{proof}
    We may assume that $Z$ is a type I or II thin set.
    When $Z$ is a type I thin set, one may argue as in \cite[Theorem 1.2]{BL19}. So we may assume that $Z = f(Y(k))$ where $f : Y  \to X$ is a dominant generically finite morphism from a smooth projective variety $Y$.
    Let $\Omega$ be a finite set of places of $k$.
    Then by Principle~\ref{conj:equidistribution}
     \begin{align*}
        & \lim_{B \to \infty}\frac{\#\{x \in f(Y(k)) \setminus \mathsf Z \, | \, \mathsf H_L(x) \leq B\} }{\#\{x \in X(k) \setminus \mathsf Z \, | \,\mathsf H_L(x) \leq B \}}\\
         &\leq \lim_{B \to \infty}\frac{\#\{x \in X(k) \setminus \mathsf Z  \, | \, x \in f(Y(k_v)) \text{ for any $v \in \Omega$}, \mathsf H_L(x) \leq B \}}{\#\{x \in X(k) \setminus \mathsf Z \, | \,\mathsf H_L(x) \leq B \}}\\
         &\leq \frac{\mathsf L_*(\Pic(X^\circ), 1)\prod_{v \in \Omega}\lambda_v^{-1}\tau_{X, E, v}(f(Y(k_v)))\prod_{v\not\in \Omega}\lambda_v^{-1}\tau_{X, E, v}(X(k_v))}{\tau_{X, E}(X^\circ (\mathbb A_k)^{\mathrm{Br}(X^\circ)})},
   \end{align*}
   where $\mathsf L_*(\Pic(X^\circ), 1)$ is the residue of the L-function of $\Pic(X^\circ)$, $\lambda_v$ is a convergence factor, and $\tau_{X, E, v}$ is a local Tamagawa measure.
   Note that we have
   \[
   \tau_{X, E, v}(f(Y(k_v))) = \tau_{X, E, v}(f(Y(\mathfrak o_v))) + O(1/q_v^{1+\epsilon})
   \]
   for some $\epsilon > 0$.
   Using \cite[Theorem 3.6.2]{Serre}, we conclude our assertion by arguing as in the proof of \cite[Theorem 1.2]{BL19}.

\end{proof}

\begin{rema}
    Without the assumption that $\mathrm{Br}(X^\circ)/\mathrm{Br}(k)$ is finite, Theorem~\ref{theo:nocontributionfromthin} fails as demonstrated by \cite[Example 8.7]{LT16}.

    Indeed, let $X$ be the singular toric Fano threefold in \cite[Example 8.7]{LT16}. Let $\beta : \widetilde{X} \to X$ be a smooth resolution and $L = -\beta^*K_X$. Then $\mathrm{Br}(\widetilde{X}^\circ)/\mathrm{Br}(k)$ is infinite and Theorem~\ref{theo:nocontributionfromthin} fails as there are infinitely many double covers $Y \to \widetilde{X}$ from a smooth projective toric threefold with the same $a, b$-invariants as $\widetilde{X}$.
    \end{rema}

\section{Construction of the surface $S$} \label{sect:surfaceconstruction}

Our construction of $S$ relies on \cite[Section 3]{DO19} and \cite{DOY23}.  The surface $S$ will admit the following features:
\begin{enumerate}
\item $S$ is a smooth projective rational surface defined over a number field $k$;
    \item $S$ admits infinitely many elliptic fibrations $\pi_i : S \to \mathbb P^1$ defined over $k$; and
    \item each $\pi_i$ admits infinitely many sections and all sections are $(-1)$-curves.
\end{enumerate}

To this end, first we construct a K3 surface $T$ with infinitely many elliptic fibrations. Then $S$ will be realized as a quotient of $T$ by an involution.

\subsection{Constructing a K3 surface}

The starting point of this construction is the following: let $k$ be a number field, and $E$ and $F$ be a pair of elliptic curves defined over $k$ with the following properties:
\begin{itemize}
    \item $E[2](k)$ and $F[2](k)$ are full; and
    \item $E$ and $F$ are geometrically non-isogenous.
\end{itemize}

\begin{rema} \label{rema:ellipticcurve}
    We note that for every number field $k$, there exists a pair of elliptic curves $E$ and $F$ satisfying the two desired properties.  Indeed, it suffices to find a suitable pair of curves over $\mathbb{Q}$.  The existence of such a pair follows from \cite[Remark 4.7]{DO19}.
\end{rema}

Let $a_1, a_2, a_3, a_4$ be the four $2$-torsion elements on $E$ and $b_1, b_2, b_3, b_4$ be the four $2$-torsion elements on $F$.
We define $T$ as the Kummer surface:
\[
T = \mathrm{Kum}(E \times F).
\]
We denote the image of $E\times\{b_i\}$ by $E^T_i$ and the image of $\{a_i\} \times F$ by $F^T_i$. Let $C^T_{i j}$ be the $(-2)$-curve lying over the nodal singularity which is the image of $(a_i, b_j)$ in the quotient surface.    Finally we let $B^{T} = \cup_{i=1}^{4} (E_{i}^{T} \cup F_{i}^{T})$.
It follows from \cite[Lemma 1.4]{Oguiso89} that for any $f^T \in \mathrm{Aut}(T)$, $f^T(B^{T}) = B^{T}$.

Following \cite[Proof of Proposition 3.7]{DO19}, we consider the following divisor:
\begin{align*}
    &D^T = E^T_1 + C^T_{11} + F^T_1 + C^T_{12} + E^T_2 + C^T_{22} + F^T_2 + C^T_{21}
\end{align*}
As explained in \cite[Proof of Proposition 3.7]{DO19}, this divisor defines an elliptic fibration:
\[
\varphi_{D^T} : T \to \mathbb P^1. 
\]
Note that $C^T_{31}, C^T_{41}$ are sections for $\varphi_{D^T}$. 
We consider $C^T_{31}$ as the zero section of $\varphi_{D^T}$. 
By \cite[Proposition 3.8]{DO19}, this fibration is relatively minimal, so that any automorphism of the generic fiber extends to the total space.
Then the translation by $C^T_{41}$ with respect to $\varphi_{D^T}$ has infinite order by \cite[After Definition 3.10]{DO19}.

Note that we can do the same construction after permuting the indices.  For any permutation $\sigma \in S_{4} \times S_{4}$, we define $D^{T,\sigma}$ to be the divisors given by the analogous sums obtained by permuting the indices by $\sigma$.

Let $\{\pi_{\ell}^T\}$ be 
\[
\{\varphi_{D^{T, \sigma}}\circ f \, |\, \sigma \in S_4 \times S_4, f \in \mathrm{Aut}(T)\}.
\]

\begin{lemm} \label{lemm:k3lemm}
    Consider the elliptic fibrations $\pi_{\ell}^{T}$.  The eight irreducible components of $B^{T}$ are the only curves that are contracted by all of these elliptic fibrations.
\end{lemm}

\begin{proof}
    Appealing to the identification of $\varphi_{D^{T}}$ as $\mathscr{J}_{1}$ in the notation of \cite[Page 669]{Oguiso89}, we see that $\varphi_{D^{T}}$ has two reducible fibers: $D^{T}$ itself and $E^T_3 + C^T_{33} + F^T_3 + C^T_{34} + E^T_4 + C^T_{44} + F^T_4 + C^T_{43}$.  By comparing these divisors against their permutations by $\sigma$, we obtain the desired statement.
\end{proof}

\subsection{Constructing the rational quotient}
Next we construct $S$. Let $\theta$ be the involution of $T$ defined by $(1_E, -1_F) \in \mathrm{Aut}(E \times F)$. 
Then by \cite[The proof of Lemma 1.4]{Oguiso89} and \cite[Lemma 2.1]{DOY23} we have
\begin{itemize}
\item for any $f^T\in \mathrm{Aut}(T)$, $f^T \theta = \theta f^T$;
    \item $T^\theta = B^{T}$; and
    \item for any smooth rational curve $R$ on $T$, $\theta(R) = R$.
\end{itemize}
Let $S$ be $T/\langle\theta\rangle$ with the projection $\rho : T \to S$. Then $S$ is a smooth rational projective surface which is the blow up of $\mathbb P^1 \times \mathbb P^1$ at $16$ points corresponding to $(a_i, b_j)$.
Let $E_i, F_i, C_{ij}, B$ be the images of $E_i^T, F_i^T, C_{ij}^T, B^{T}$ respectively.
We define the divisor
\begin{align*}
    &D = E_1 + 2C_{11} + F_1 + 2C_{12} + E_2 + 2C_{22} + F_2 + 2C_{21}
\end{align*}
Then we have $\rho^*D = 2D^{T}$ so that this divisor defines an elliptic fibration:
\[
\varphi_{D} : S \to \mathbb P^1. 
\]
Again $C_{31}, C_{41}$ are sections for $\varphi_{D}$.  
We consider $C_{31}$ as the zero section of $\varphi_{D}$.
By descent from $T$ we see that the automorphisms of the generic fiber of $\varphi_D$ act as automorphisms on $S$.
Just as before, the translation by $C_{41}$ with respect to $\varphi_{D}$ has infinite order. 
All the statements still hold true after permuting the indices by $\sigma \in S_{4} \times S_{4}$: we get a divisor $D^{\sigma}$ and an elliptic fibration $\varphi_{D^{\sigma}}$.

Let $\{\pi_{\ell}\}$ be 
\[
\{\varphi_{D^\sigma}\circ f \, |\, \sigma \in S_4 \times S_4, f \in \mathrm{Aut}(S)\}.
\]
Henceforth we denote the infinite set of sections of all the elliptic fibrations $\pi_{\ell}$ by $\{ C_{i} \}$. Each $C_{i}$ is a $(-1)$-curve.

\begin{lemm} \label{lemm:infiniteorbit}
    Suppose that $C$ is an irreducible component of a reducible fiber of $\pi_{\ell}$ for some $\ell$.
    \begin{enumerate}
        \item If $C \subset B$, then every element in $\Aut(S)$ sends $C$ to another curve in $B$. 
        \item If $C \not \subset B$, then the orbit of $C$ under $\Aut(S)$ is infinite.
    \end{enumerate}
\end{lemm}

\begin{proof}
    (1) Using the blow-up description of $S$ it is easy to see that $\sum_i (E_i + F_i)$ is contained in $|-2K_{S}|$.  Since $-2K_{S}$ intersects each curve $E_{i},F_{i}$ negatively, in fact $\sum_i (E_i + F_i)$ is the unique curve in this linear series and the statement follows. 

    (2) Since each of the fibrations $\pi_{\ell}$ can be transformed to a fibration $\varphi_{D^{\sigma}}$ via an automorphism of $S$, it suffices to prove the statement for the fibration $\varphi_{D^{\sigma}}$. The reducible fibers (of the corresponding map on $T$) are described by \cite[Page 669]{Oguiso89}; each irreducible component either lies in $B$ or has the form $C_{ij}$ for some $i,j$.  But each $C_{ij}$ is a section of a fibration $\varphi_{D^{\widehat{\sigma}}}$ for a different choice of permutation $\widehat{\sigma}$.  Thus it will have infinite orbit under the subgroup generated by the automorphisms of the generic fiber of $\varphi_{D^{\widehat{\sigma}}}$.
\end{proof}

\begin{prop} \label{prop:coversofcurves}
Suppose that $\{ g_{j}: W_{j} \to S\}$ is a finite set of dominant finite thin maps such that every $W_{j}$ is normal and geometrically irreducible.  Then for infinitely many $i$, the preimage of $C_{i}$ in every $W_{j}$ is a union of curves which each have degree $\geq 2$ over $C_{i}$.
\end{prop}

\begin{proof}
    Let $B_{j} \subset S$ denote the branch divisor of $g_{j}$.  Since $S$ is rational, every $B_{j}$ must be non-zero.  We separate into two cases:
    
    \textbf{Case 1:} $\Supp(B_{j})$ is contained in $B$ for some $j$.
    
    The geometric fundamental group of $T \backslash B^{T}$ is $\mathbb{Z}/2\mathbb{Z}$ by \cite[Table 1, Number 3]{KZ02}.  We claim that the geometric fundamental group of $S \backslash B$ is isomorphic to $\mathbb{Z}/4\mathbb{Z}$.  Indeed, $S \backslash B$ contains the open set $V$ which is the complement of $8$ vertical fibers $\mathbb P^1\times\{b_i\}$ and $\{a_i\}\times \mathbb P^1$ in $\mathbb{P}^{1} \times \mathbb{P}^{1}$.  The geometric fundamental group of $V$ is the product
    \begin{equation*}
        \langle \alpha_{1},\alpha_{2},\alpha_{3},\alpha_{4} | \alpha_{1}\alpha_{2}\alpha_{3}\alpha_{4}=1 \rangle \times \langle \beta_{1},\beta_{2},\beta_{3},\beta_{4} | \beta_{1}\beta_{2}\beta_{3}\beta_{4}=1 \rangle
    \end{equation*}
    where $\alpha_{i}$ is a loop around $a_{i}$ and $\beta_{j}$ is a loop around $b_{j}$.  By adding back in the curves $C_{ij} \backslash B$, we kill all the products $\alpha_{i}\beta_{j}$.  Thus $\pi_{1}((S\setminus B)_{\overline{k}})$ is generated by a single class $\alpha$ satisfying $\alpha^{4}=1$.

    In particular, $S_{\bar{k}}$ admits two different covers $T_{\bar{k}}, U_{\bar{k}}$ that are ramified only over $B$.  After base-change to the algebraic closure, the map $g_{j}$ must be one of these finite covers.
    
    Recall that each map $C_{i}^{T} \to C_{i}$ is a double cover.  Since $U_{\bar{k}} \to S_{\bar{k}}$ factors through $T_{\bar{k}}$, it is also true that every irreducible component of the preimage of $C_{i,\bar{k}}$ in $U_{\bar{k}}$ has degree $\geq 2$ over $C_{i,\bar{k}}$.  This property descends to $W_{j}$.  Thus $W_{j}$ has the desired property for every curve $C_{i}$.  We can repeatedly remove such $W_{j}$ from our set of thin maps so that we land in Case 2.
    
    \textbf{Case 2:} $\Supp(B_{j})$ is not contained in $B$ for any $j$.

    Let $\{ L_{k} \}$ denote the set of irreducible components of the various $B_{j}$ that are not contained in $B$.  Note that there are only finitely many such curves.
    
    We claim that there is an elliptic fibration $\pi: S \to \mathbb{P}^{1}$ such that none of the $L_{k}$ are contracted by $\pi$.  
    Indeed, fix any elliptic fibration $\pi_{1}: S \to \mathbb{P}^{1}$ and let $\{Q_\ell\}$ denote the irreducible components of fibers that are not contained in $B$. (For irreducible fibers, we represent all of such components by the class $Q_\ell$ of irreducible fibers.)  For each such curve $Q_\ell$, the stabilizer $G_{\ell}$ is a subgroup of $\Aut(S)$ with infinite index.  For each curve $L_{k}$, let $G^{k}_{\ell}$ denote the coset of $G_{\ell}$ that sends $Q_{\ell}$ to $L_k$ when such a coset exists.  Due to the infinite index property, \cite[Lemma 4.1]{Neumann54} shows that the union $\cup_{\ell,k} G^{k}_{\ell}$ is a proper subset of $\Aut(S)$.  Thus, by acting on $\pi_{1}$ by elements of $\Aut(S)$ not contained in this union of cosets we obtain infinitely many elliptic fibrations that do not contain any curve $L_{k}$ as an irreducible component in a reducible fiber.  Some of these infinitely many such fibrations will also have the property that no curve $L_{k}$ occurs as an irreducible fiber, yielding the desired map $\pi$.  Since each branch divisor $B_{j}$ contains at least one of the $L_{k}$, we see that $B_{j}$ has positive intersection against the general fiber of $\pi$. 

    We claim that there are infinitely many $C_{i}$ which are sections of $\pi$ which have the desired property.  First suppose that the composition $\pi \circ g_{j}: W_{j} \to \mathbb{P}^{1}$ fails to have connected fibers.  Let $\psi_{j}: V_{j} \to \mathbb{P}^{1}$ denote the finite part of the Stein factorization.  Then any curve mapping to a section $C_{i}$ of $\pi$ must have degree at least $\deg(\psi_{j})$.

    Next suppose that $\pi \circ g_{j}$ has connected fibers.  
    Let $E$ denote the generic fiber of $\pi$ and let $R_{j}$ denote the generic fiber of $\pi \circ g_j$, both defined over the field $K = K(\mathbb{P}^{1}_{k})$.  Since $W_{j}$ is geometrically irreducible, $R_{j}$ is geometrically integral.  Furthermore, by construction the map $R_{j} \to E$ has a non-trivial branch divisor.  By \cite[Theorem 1.3]{CDJLZ22} there is a finite index coset $Q \subset E(K)$ such that for every $j$ the intersection $Q \cap g_{j}(R_{j}(K))$ is empty.  Since every element of $Q$ defines a $(-1)$-curve, this proves the desired statement for the infinitely many $(-1)$-curves in the coset $Q$.
\end{proof}

\section{Construction of the threefold $X$}

We next discuss the construction of the threefold $X$.  Choose an ample divisor $H$ on $S$ such that $H - K_{S}$ is big.  We define
\begin{equation*}
    X = \mathbb{P}_{S}(\mathcal{O}_{S} \oplus \mathcal{O}_{S}(-H)).
\end{equation*}
We denote the $\mathbb{P}^{1}$-bundle structure by $\pi: X \to S$.  The anticanonical divisor is $-K_{X} = 2 \xi + \pi^{*}(H-K_{S})$ where $\xi$ denotes a divisor representing $\mathcal{O}_{X/S}(1)$.  By \cite[Lemma 2.3.2]{Lazarsfeld04} $-K_{X}$ is big.

Since $S$ is rational, so is $X$.  This implies that $\rH^{1}(X,\mathcal{O}_{X}) = 0$, $\rH^{2}(X,\mathcal{O}_{X}) = 0$, and that $\Pic(X_{\overline{k}})$ is torsion-free, since all these properties are birational invariants for smooth projective varieties.  Thus $X$ is almost Fano.

\subsection{Counting points}

We let $T_{i} \subset X$ denote the preimage of $C_{i} \subset S$ as defined in Section \ref{sect:surfaceconstruction}.

\begin{lemm} \label{lemm:thinsubsetintersects}
    Suppose that $Z \subset X(k)$ is a thin set.  Then there are infinitely many $i$ such that the set $Z \cap T_{i}(k)$ is a thin subset of $T_{i}(k)$.
\end{lemm}

\begin{proof}
    It suffices to consider the contributions of a finite set of dominant thin maps $\{f_{j}: Y_{j} \to X\}$ where $Y_{j}$ is normal and geometrically irreducible.  Define the composition $h_{j} = \pi \circ f_{j}: Y_{j} \to S$. We let $g_{j}: W_{j} \to S$ denote the finite part of the Stein factorization of $h_{j}$.  We apply Proposition \ref{prop:coversofcurves} to the finite set of maps $\{g_{j}: W_{j} \to S | \deg(g_{j}) > 1 \}$ to obtain an infinite set of curves $C_{i} \subset S$.  We let $T_{i} \subset X$ denote the preimage of $C_{i}$. 

    First consider an index $j$ such that $h_{j}$ fails to have connected fibers.  By Proposition \ref{prop:coversofcurves}, the preimage of $C_{i}$ in $W_{j}$ is a union of curves which map to $C_{i}$ with degree $\geq 2$.  Thus the preimage of $T_{i}$ in $Y_{j}$ consists of a union of surfaces of degree $\geq 2$ over $T_{i}$.  We conclude that $f_{j}(Y_{j}(k)) \cap T_{i}(k)$ is thin in $T_{i}(k)$.

    Next consider an index $j$ such that $h_{j}$ has connected fibers.  For every $s$ outside of a closed subset of $S$, the fiber of $h_{j}$ over $s$ is a geometrically irreducible curve.  Thus, for all but finitely many $i$ the preimage of $C_{i}$ in $Y_{j}$ is a geometrically irreducible surface $Q_{i,j}$ that is not contained in the ramification locus.  This means that $Q_{i,j} \to T_{i}$ is a thin map, showing that $f_{j}(Y_{j}(k)) \cap T_{i}(k) = f_{j}(Q_{i,j}(k))$ is a thin subset of $T_{i}(k)$.
\end{proof}

\begin{proof}[Proof of Theorem \ref{theo:maintheorem}:]
   Each $T_{i}$ is isomorphic to a Hirzebruch surface $\mathbb{F}_{e_{i}} := \mathbb{P}_{\mathbb{P}^{1}}(\mathcal{O}_{\mathbb{P}^{1}} \oplus \mathcal{O}_{\mathbb{P}^{1}}(-e_{i}))$ where $e_{i} = H \cdot C_{i}$.
    Note that $K_{T_{i}} = -2\xi|_{T_i} + (-e_{i}-2)\pi^{*}\mathcal{O}_{C_{i}}(1)$.  Thus
    \begin{equation*}
        K_{T_{i}} - \tfrac{e_{i}+2}{e_{i}+1}  K_{X}|_{T_i}  = \tfrac{2}{e_{i}+1}\xi|_{T_i}
    \end{equation*}
    lies on the boundary of $\Eff^{1}(T_{i})$ showing that $a(T_{i},-K_{X}) = \frac{e_{i}+2}{e_{i}+1}$. 
    Since $T_{i}$ is toric, for any thin set $Z_{i} \subset T_{i}(k)$ Theorem~\ref{theo:santens}, Theorem~\ref{theo:nocontributionfromthin}, and \cite{BT96b,BT98,CLT01} 
    show that there is some constant $\mathsf C_i > 0$ such that
    \begin{align*}
        \# \{ \, x \in T_{i}(k) \backslash Z_{i} \,\,| \,\, \mathsf H_{-K_{X}}(x) \leq B \, \} \geq \mathsf C_{i}B^{\frac{e_{i}+2}{e_{i}+1}}(1   + o(1)).
    \end{align*}

    Suppose that $Z \subset X(k)$ is any thin set.  By Lemma \ref{lemm:thinsubsetintersects}, we see that there is some $T_{i}$ such that $Z \cap T_{i}(k)$ is a thin subset of $T_{i}(k)$.  If we set $\mathsf C = \tfrac{\mathsf C_{i}}{2}$ and $\epsilon = \tfrac{1}{e_{i}+1}$ then we conclude that
        \begin{align*}
        \#\{ \, x \in X(k) \backslash Z \, \,  | \, \, \mathsf H_{-K_{X}}(x) \leq B \, \} & \geq \# \{ \, x \in T_{i}(k) \backslash Z \,\,| \,\, \mathsf H_{-K_{X}}(x) \leq B \, \} \\
        & > \mathsf CB^{1+\epsilon}
    \end{align*}
    for all sufficiently large $B$.
\end{proof}

\bibliographystyle{alpha}
\bibliography{almostfanoex}

\end{document}